\documentclass[a4paper,12pt]{article}
\usepackage{amssymb}
\usepackage{epsfig}
\usepackage{amsfonts}
\usepackage{amsmath}
\usepackage{euscript}
\usepackage{amscd}
\usepackage{amsthm}
\usepackage{hyperref}
\usepackage{xcolor}
\DeclareMathAlphabet{\mathpzc}{OT1}{pzc}{m}{it}

\newtheorem{thm}{Theorem}[section]
\newtheorem{lem}[thm]{Lemma}

\newtheorem{prop}[thm]{Proposition} 

\newtheorem{defn}[thm]{Definition}

\newtheorem{cor}[thm]{Corollary}

\newtheorem{rem}[thm]{Remark}

\newtheorem{ex}[thm]{Example}

\newcommand{\bZ}{\mathbb Z}
\newcommand{\bQ}{\mathbb Q}
\newcommand{\bN}{\mathbb N}
\newcommand{\bC}{\mathbb C}
\newcommand{\bK}{\mathbb K}
\newcommand{\Ga}{\mathbb G}
\newcommand{\V}{\mathbb V}
\newcommand{\X}{\mathbb X}
\newcommand{\bD}{\mathbb D}
\newcommand{\bT}{\mathbb T}

\newcommand{\U}{\mathbb U}
\newcommand{\bS}{\mathbb S}
\newcommand{\bH}{\mathbb H}

\newcommand{\rbo}{\big{(}}
\newcommand{\rbc}{\big{)}}
\author{Nikhilesh Dasgupta$^{*}$ and Animesh Lahiri$^\dagger$\\
{\small{\it $^{*}$ National Research University Higher School of Economics,}}\\
{\small{\it Faculty of Computer Science,}}\\
{\small{\it Pokrovsky Bulvar 11, Moscow 109208, Russia.}}\\
{\small{\it e-mail: its.nikhilesh@gmail.com, ndasgupta@hse.ru}}\\
{\small {\it $^\dagger$ Chennai Mathematical Institute,}}\\
{\small {\it H1 SIPCOT IT Park, Siruseri, Kelambakkam, Chennai 603103, India. }}\\
{\small {\it e-mail : animeshl@cmi.ac.in}}}
\title{Isotropy subgroups of some almost rigid domains}

\begin{document}
\date{} 
\maketitle
\abstract{ 

\noindent 
In this paper, we describe the structure of isotropy subgroups of some almost rigid domains and give necessary and sufficient conditions for an automorphism to be an element of the isotropy subgroup.

\smallskip
\noindent
{\small {{\bf Keywords}. Locally nilpotent derivation, almost rigid domains, automorphism group, isotropy subgroup.}

}}

\section{Introduction}
Let $k$ be a field of characteristic zero and $B$ a $k$-algebra. A $k$-linear map $\delta:B\longrightarrow B$ is said to be a {\it $k$-derivation}, if it satisfies the Leibnitz rule: $\delta(ab)=a\delta(b)+b\delta(a),~\forall a,b \in B$. A $k$-derivation $\delta$ is said to be a {\it locally nilpotent derivation} (abbrev. {\it lnd}) if, for each $b\in{B}$, there exists $n(b)\in{\bN}$ such that $\delta^{n(b)}(b)=0$. The set of all locally nilpotent $k$-derivations on $B$ will be denoted by ${\rm LND}(B)$ and the set of all $k$-algebra automorphisms of $B$ will be denoted by ${\rm Aut}(B)$. The {\it kernel} of $\delta$, denoted by ${\rm Ker}(\delta)$, is defined to be the set $\{b \in B ~|~\delta(b)=0\}$.

\medskip
\indent
There is a natural action of ${\rm Aut}(B)$ on ${\rm LND}(B)$ defined by $\alpha\cdot \delta=\alpha\delta{\alpha}^{-1}$, for $\alpha\in{\rm Aut}(B)$ and $\delta\in{\rm LND}(B)$. For a given $\delta\in{\rm LND}(B)$, the {\it stabilizer of $\delta$} under the above action, i.e., the subgroup 
$\lbrace{\sigma}\in{\rm Aut}(B):\sigma \delta=\delta\sigma\rbrace$ of ${\rm Aut}(B)$, is called the {\it isotropy subgroup} of $B$ with respect to $\delta$ and will be denoted by ${\rm Aut}(\delta)$. Every $\delta \in {\rm LND}(B)$ induces an element of ${\rm Aut}(B)$ via the exponential map, defined as ${\rm Exp}(\delta):=\underset{i \geqslant 0}{\sum}\dfrac{1}{i!}{\delta}^i$. For any {\it lnd} $\delta$, each of its {\it replicas} $f\delta$ \big{(}$f \in {\rm Ker}(\delta)$\big{)}, is also an {\it lnd}. Exponents of all replicas of $\delta$ form a commutative subgroup $\U(\delta):=\{{\rm Exp}(f\delta)~|~f \in {\rm Ker}(\delta)\}$, called the {\it big unipotent subgroup corresponding to $\delta$}. The correspondence $f \leftrightarrow {\rm Exp}(f\delta)$ induces an isomorphism between $\U(\delta)$ and $\big{(}{\rm Ker}(\delta), +\big{)}$. It is easy to see that for any $\delta \in {\rm LND}(B)$, the big unipotent group $\U(\delta)$ is a subgroup of ${\rm Aut}(\delta)$. This information motivates one to investigate when for a $k$-algebra $B$ and an {\it lnd} $\delta$, ${\rm Aut}(\delta)$ properly contains $\U(\delta)$.

\smallskip
\noindent
 Clearly, $\U(\delta)$ is properly contained in ${\rm Aut}(\delta)$ when $B$ admits an {\it lnd} ${\delta}^{\prime}$ which commutes with $\delta$ and ${\rm Ker}(\delta) \neq {\rm Ker}({\delta}^{\prime})$ (for example, when $B=k[X,Y]$ and $\delta=\frac{\partial}{\partial X}$). Indeed, if such a ${\delta}^{\prime}$ exists, then, for any $f \in {\rm Ker}(\delta) \cap {\rm Ker}({\delta}^{\prime})$, ${\rm Exp}(f{\delta}^{\prime})$ is an element of ${\rm Aut}(\delta)\setminus\U(\delta)$ (see \hyperref[al2]{\color{blue}{Remark \ref{al2}}}). So, it is interesting to investigate affine domains where all non-zero locally nilpotent derivations have the same kernel.
 
\medskip
\indent
The goal of this paper is to study the isotropy subgroups of locally nilpotent derivations on some {\it almost rigid domains}. An affine $k$-domain $B$ is said to be {\it almost rigid} if there exists $D \in {\rm LND}(B)$ such that every $\delta \in {\rm LND}(B)$ can be written as $\delta =hD$, for some $h \in {\rm Ker}(D)$. Moreover, $D$ is called a {\it canonical lnd} on $B$. For an almost rigid domain $B$ with no non-constant invertible elements, if $D$ is a canonical {\it lnd} and $\phi \in {\rm Aut}(B)$, then $\phi D {\phi}^{-1}=\lambda D$, for some $\lambda \in k^*$ (see \hyperref[al1]{\color{blue}{Lemma \ref{al1}}}). This motivates one to investigate which automorphisms of $B$ are in ${\rm Aut}(D)$.

\smallskip
\noindent
Recently, R. Baltazar and M. Veloso (\cite{BV}) have studied the isotropy subgroups of  certain types of Danielewski surfaces, which are examples of some almost rigid domains. In fact, they have described the isotropy subgroups of locally nilpotent derivations on Danielewski surfaces defined by polynomials of the form $f(X)Y-P(Z)$ over an algebraically closed field of characteristic zero. In this paper, we have studied the structure of the isotropy subgroup of an arbitrary {\it lnd} on more general surfaces and also on almost rigid domains of Krull dimension more than two.

\medskip
\indent
In Section \hyperref[iso_GDs]{\color{blue}\ref{iso_GDs}}, we have studied the isotropy subgroups of generalized Danielewski surfaces (introduced by A. Dubouloz in \cite{D}). For any {\it lnd}, we have shown that the isotropy subgroup is isomorphic to a semi-direct product of the big unipotent group, corresponding to a {\it canonical lnd}, with a commutative subgroup containing automorphisms corresponding to a non-trivial action of a one-dimensional algebraic torus (see \hyperref[thmdp]{\color{blue}{Theorem \ref{thmdp}}}). We have also given necessary and sufficient conditions on an automorphism to be an element of the isotropy group with respect to {\it any} replica of the {\it canonical lnd} (see \hyperref[y_1DS]{\color{blue}Lemma \ref{y_1DS}}).

\smallskip
\noindent
 In Sections \hyperref[iso_Dvc]{\color{blue}\ref{iso_Dvc}} and \hyperref[iso_Dv]{\color{blue}\ref{iso_Dv}}, Danielewski varieties (introduced by A. Dubouloz in \cite{D}), which are further generalizations of Danielewski surfaces to higher dimensions, have been considered and their isotropy subgroups have been studied. For Danielewski varieties with {\it constant coefficients}, the isotropy subgroup of a {\it canonical lnd} is isomorphic to a semi-direct product of a normal subgroup containing the big unipotent group and automorphisms induced by an algebraic torus, with a group of permutations (see \hyperref[danst2]{\color{blue}Theorem \ref{danst2}}). For {\it general} Danielewski varieties, we have proved that the isotropy subgroup of the {\it canonical lnd} is isomorphic to a semi-direct product of the big unipotent group with a suitable subgroup of the automorphism group (see  \hyperref[danv3]{\color{blue}Corollary \ref{danv3}}). Moreover, in each of these cases, we have given necessary and sufficient conditions on an automorphism to be an element of the isotropy group with respect to {\it any} replica of the canonical {\it lnd}.

\smallskip
\noindent 
In case of generalized Danielewski surfaces (see \hyperref[dansaut]{\color{blue}Remarks \ref {dansaut}}, \hyperref[dansrem] {\color{blue}\ref{dansrem}}) and Danielewski varieties with {\it constant coefficients}  (see \hyperref[ptorus2]{\color{blue}Corollary \ref{ptorus2}}, \hyperref[dancrem]{\color{blue}Remark \ref{dancrem}}), the isotropy subgroup ${\rm Aut}(\delta)$ always contains automorphisms induced by actions of an algebraic torus, when $k$ is algebraically closed and $\delta$ is the canonical {\it lnd}. So, in these cases, the isotropy subgroup properly contains the big unipotent group. However, there are examples of Danielewski varieties (with non-constant coefficients), where the isotropy group is equal to the unipotent group $\U(\delta)$ (see \hyperref[danrem]{\color{blue}Remark \ref{danrem}}). 
  
\smallskip
\noindent  
In Section \hyperref[iso_FM]{\color{blue}\ref{iso_FM}}, we have considered a threefold over $\bC$ given by D.R. Finston and S.~ Maubach \cite{FM} and have described the structure of ${\rm Aut}(\delta)$ for any locally nilpotent derivation $\delta$ (see \hyperref[thmfm]{\color{blue}Theorem \ref{thmfm}}). Finally, in Section \hyperref[iso_DDs]{\color{blue}\ref{iso_DDs}}, isotropy subgroups of Double Danielewski surfaces (introduced by N Gupta and S. Sen in \cite{GS1}) have been studied and we have given necessary and sufficient conditions for an automorphism to be an element of the isotropy subgroup with respect to any replica of the canonical {\it lnd} (see \hyperref[y_1]{\color{blue}Proposition \ref{y_1}}).   

\section{Notations}
By a ring $R$, we will mean a commutative ring with unity. Unless stated otherwise, capital letters like $Y_1,Y_2,\dots,Y_m,$ $X,Y,Z$ will be used as indeterminates over respective ground rings. Throughout the paper, $k$ (respectively $\bK$) will stand for a field (respectively algebraically closed) of characteristic zero. For a ring $R$, $R^*$ denotes the set of all units of $R$. For $f\in{R[X_1,\dots,X_n]}$, we denote the {\it support of $f$} by ${\rm Supp}(f)$. When we say $\alpha=[F_1,\dots,F_n]$ is an automorphism of $R[X_1,\dots,X_n]$, we will mean that $\alpha$ is an automorphism such that $\alpha(X_i)=F_i$, for $1\leqslant i\leqslant n$.

\smallskip
\noindent
Let $B$ be a $k$-domain. For any $\delta \in {\rm LND}(B)$ and a subgroup $G$ of ${\rm Aut}(B)$, let us define 
$$G_{\delta}:=G \cap {\rm Aut}(\delta).$$ 
The Makar-Limanov invariant of $B$, denoted by ${\rm ML}(B)$, is defined as follows :
$${\rm ML}(B)=\bigcap_{\delta \in {\rm LND}(B)} {\rm Ker}(\delta).$$ It is easy to see that if $B$ is an almost rigid domain with canonical {\it lnd} $\delta$, then ${\rm ML}(B)={\rm Ker}(\delta).$

\smallskip
\noindent
For a group $G$, a subgroup $H$ of $G$ and an element $g \in G$, we use the notations $C_G(H)$ and $C_G(g)$ to denote the centralizers (in $G$) of $H$ and $g$ respectively.    

\medskip
\indent
The rest of the paper has been divided into six sections. We begin by proving some basic results on isotropy subgroups and almost rigid domains in the first section. In the subsequent sections, we have studied the isotropy subgroups of locally nilpotent derivations on different classes of almost rigid domains.

\section{Some basic results}\label{basic}
\begin{prop}\label{iso1} Let $B$ be a $k$-domain, $\delta (\neq 0) \in {\rm LND}(B)$ and $A={\rm Ker}(\delta)$. For any $f \in A$, let $H_f:=\{\theta \in {\rm Aut}(\delta)~|~\theta(f)=f\}$. Then the following statements hold.
\begin{enumerate}
\item [\rm (i)] For any $\theta \in {\rm Aut}(B)$, if the lnd $\theta \delta {\theta}^{-1} \in  \U(\delta)$, then $A$ is $\theta$-invariant. Moreover, ${\theta\vert_{A}}\in {\rm Aut}(A)$. 
\item [\rm (ii)] $C_{{\rm Aut}(\delta)}\big{(}{\rm Exp}(f\delta)\big{)}=H_f$.
\item [\rm (iii)] $C_{{\rm Aut}(\delta)}\rbo\U(\delta)\rbc=\{\theta \in {\rm Aut}(\delta)~|~\theta\vert_A={Id}_A\}.$ 
\end{enumerate}
\end{prop}

\begin{proof} (i) Let $\theta \delta {\theta}^{-1}=g\delta$, for some $g(\neq 0) \in A$. Let $a \in A$. Since $\theta$ is onto, there exists $b \in B$, such that $\theta (b)=a$. Then 
$$0=g\delta(a)=\theta \delta {\theta}^{-1}(a)=\theta (\delta (b)).$$ Since $\theta \in {\rm Aut}(A)$, we have $\delta (b)=0$, i.e., $b \in A$. Hence $A \subseteq \theta(A)$. Conversely, for any $a \in A$, 
$$g\delta (\theta (a))=\theta \delta (a)=0.$$ So $\delta (\theta (a))=0$ and hence $\theta(A) \subseteq A$.

\smallskip
\noindent
(ii)  Suppose $\theta\in{H_f}$. Then, for any $i \in \bN$, $\theta {\delta}^i={\delta}^i \theta $  and $\theta (f^i)=f^i$. Now, for any $i \in \bN$ and $x \in B$, we have
$$\theta (f\delta)^i(x)=\theta\rbo f^i{\delta}^i(x)\rbc=f^i\theta {\delta}^i(x)=f^i{\delta}^i\theta(x)=(f\delta)^i\theta(x)$$
and hence $\theta$ commutes with ${\rm Exp}(f\delta)$.

\smallskip
\noindent
 Conversely, suppose $\theta\in{\rm Aut}(\delta)$ and $\theta {\rm Exp}(f\delta) ={\rm Exp}(f\delta) \theta$. Since $\delta \neq 0$, there exists $s \in B$ such that $\delta (s) \neq 0$ and ${\delta}^2(s)=0$. Then 
$$\theta {\rm Exp}(f\delta)(s)=\theta (s+f\delta(s))=\theta(s)+\theta(f)\theta\delta(s),$$ and
$${\rm Exp}(f\delta)\theta (s)=\theta(s)+f\theta\delta (s).$$ Since $B$ is a domain, we have 
$\theta (f)=f$.

\smallskip
\noindent
(iii) The proof follows from the following observation :
$$C_{{\rm Aut}(\delta)}\big{(}\U(\delta)\big{)}=\bigcap_{f \in A} C_{{\rm Aut}(\delta)}\big{(}{\rm Exp}(f \delta)\big{)}=\bigcap_{f \in A} H_f=\{ \theta \in {\rm Aut}(\delta)~|~\theta|_A={Id}_A \}.$$ 

\end{proof}

\begin{rem}\label{genrem} {\em For an almost rigid domain $B$ with canonical {\it lnd} $\delta$ and $A={\rm Ker}(\delta)$, if $\phi \in {\rm Aut}(B)$, then $\phi\vert_A \in {\rm Aut}(A)$ by Proposition \ref{iso1} (i).
}
\end{rem}

\begin{lem}\label{equilnd}
Let $B$ be a $k$-domain and $D_1,D_2\in{\rm LND}(B)$ and $f\in {\rm Ker}(D_1) \cap{\rm Ker}(D_2)$. If $D_1D_2=D_2D_1$, then ${\rm Exp}(fD_1)\in{\rm Aut}(D_2).$ In particular, the result holds if ${\rm Ker}(D_1)={\rm Ker}(D_2)$.
\end{lem}
\begin{proof}  Let $x \in B$. Then for any $i \in \bN $, 
$$(fD_1)^iD_2(x)=f^i(D_1^iD_2)(x)=f^i(D_2D_1^i)(x)=D_2(f^iD_1^i)(x)=D_2(fD_1)^i(x),$$
where the second equality holds because $D_1D_2=D_2D_1$. Hence ${\rm Exp}(fD_1)\in{\rm Aut}(D_2).$ If ${\rm Ker}(D_1)={\rm Ker}(D_2)(=A)$, then By Principle 12 in \cite{F}, there exists $a,b\in A$ such that $aD_1=bD_2$. So $D_1D_2=D_2D_1$ and hence the result follows.
\end{proof}

\begin{rem}\label{al2}{\em From Lemma \ref{equilnd}, we note that if $D_1$ and $D_2$ are two commuting locally nilpotent derivations with ${\rm Ker}(D_1) \neq {\rm Ker}(D_2)$, then for any $f \in {\rm Ker}(D_1) \cap{\rm Ker}(D_2)$, ${\rm Exp}(fD_1)$ is an element of ${\rm Aut}(D_2)\setminus\U(D_2)$.
}
\end{rem}

\indent
The next result describes the orbit of the action of the automorphism group on the canonical {\it lnd} for an almost rigid domain.
\begin{lem}\label{al1}
Let $B$ be an almost rigid domain with canonical {\it lnd} $D$, $A={\rm Ker}(D)$ and $\phi \in {\rm Aut}(B)$. Then, $\phi D {\phi}^{-1} =fD$, for some $f \in A^*$.
\end{lem}

\begin{proof} Since $B$ is almost rigid, there exist $f\in A$ such that $\phi D {\phi}^{-1}=fD$. Then 
$$D={\phi}^{-1}(fD)\phi={\phi}^{-1}(f)({\phi}^{-1}D \phi)={\phi}^{-1}(f)gD,~~\text{for some}~g \in A.$$
Since $B$ is a domain and $A$ is factorially closed in $B$, $f \in A^*$.
\end{proof}

\indent
Now, we prove two elementary technical lemmas which will be used throughout rest of the paper.

\begin{lem}\label{root of unity}
Let $f\in k[X_1,\dots,X_n]$. If $\lambda ,\mu\in{k^*}$ be such that $\lambda f(X_1, \dots ,X_l)={\mu}^mf({\mu}^{a_1} X_1,\dots,{\mu}^{a_l}X_l)$, for some $m,a_1,\dots,a_l\in{\bZ}$, then ${\mu}^{m+a_1i_1+\dots+a_li_l}=\lambda$, for each $(i_1,\dots,i_l)\in{\rm Supp}(f)$.
\end{lem}
\begin{proof}
 Let $f(X_1,\dots ,X_l)={\underset{(i_1,\dots ,i_l)\in{\rm Supp}(f)}{\sum}}{a_{i_1,\dots,i_l}X_1^{i_1}\dots X_l^{i_l}}\in{k[X_1, \dots,X_l]}$. Since $\lambda \Big{(}{\underset{(i_1,\dots ,i_l) \in{\rm Supp}(f)}{\sum}}{a_{i_1,\dots,i_l}X_1^{i_1}\dots X_l^{i_l}}\Big{)}={\underset{(i_1,\dots ,i_l) \in{\rm Supp}(f)}{\sum}}{{\mu}^{m+a_1i_1+\dots+a_li_l}a_{i_1,\dots,i_l}X_1^{i_1}\dots X_l^{i_l}}$, we have ${\mu}^{m+a_1i_1+\dots+a_li_l}=\lambda$, for each $(i_1,\dots,i_l)\in{\rm Supp}(f)$.
\end{proof}

\begin{lem}{\label{chainrule}}
Let $R$ be a ring containing $\bQ$, $P \in B=R[X_1,\dots,X_n]$, $\delta\in{\rm LND}(B)$ and $\alpha\in{\rm Aut}_{R}(B)$. Then $\delta\Big{(}P\big{(}\alpha(X_1),\dots,\alpha(X_n)\big{)}\Big{)}=\overset{n}{\underset{i=1}\sum}\alpha\Big{(}\dfrac{\partial P}{\partial {X_i}}\Big{)}\delta\big{(}\alpha(X_i)\big{)}$.
\end{lem}

\begin{proof} The proof follows from the fact that since $\alpha \in {\rm Aut}(B)$, for each $i \in \{1, \dots ,n\}$, we have $\dfrac{\partial P(\alpha(X_1), \dots ,\alpha(X_n))}{\partial \alpha(X_i)}=\alpha\Big{(}\dfrac{\partial P}{\partial X_i}\Big{)}$.
\end{proof}

\section{Generalised Danielewski surfaces}\label{iso_GDs}
In this section we consider the generalised Danielewski surfaces introduced by A.~ Dubouloz in \cite{D}. We describe the isotropy subgroup of an arbitrary {\it lnd} on a Danielewski surface which is in a standard form. Let us first recall some definitions and fix some notations which will be used throughout this section.  

\begin{defn}\label{Ds}
{\em Let $B$ be a $k$-algebra. $B$ is said to be a {\it generalised Danielewski surface over $k$} if $B$ is isomorphic to the $k$-algebra $$B_{d,P}:=\dfrac{k[X,Y_1,Y_2]}{\big{(}X^dY_2-P(X,Y_1)\big{)}},$$ where $d\geqslant{2}$ and $r:={\rm deg}_{Y_1}(P)\geqslant{2}$. If $P(X,Y_1)={\underset{i=1}{\overset{r}\prod} \big{(}Y_1-\sigma_i(X)\big{)}}$, where $\sigma_i(X)\in k[X]$, then the surface $B_{d,P}$ is called a generalised Danielewski surface in {\it standard form}.
}
\end{defn}

\smallskip
\noindent
 By $x,y_1,y_2$ we will denote the images of $X,Y_1,Y_2$, respectively, in $B_{d,P}$.
%
The automorphism group of $B_{d,P}$ was studied by A. Dubouloz and P-M. Poloni (\cite[Theorem 3.11 and Lemma 3.14]{DP}). Let $\mathcal{S}_r$ denote the symmetric group of $r$ elements and $id$ denotes the identity permutation.

\smallskip
\noindent
 Every automorphism $\Phi$ in ${\rm Aut}(B_{d,P})$ is uniquely determined by the datum $\mathcal{A}_{\Phi}=\big{(}\alpha,\mu,a, b(x)\big{)}\in {\mathcal{S}_r\times k^*\times k^*\times k[x]}$, such that the polynomial $c(x):=\sigma_{\alpha(i)}(ax)-\mu\sigma_i(x)$ does not depend on the index $i=1,2,\dots,r$. Moreover, ${\Phi}$ is induced by the automorphism $\Psi$ of $k[X,Y_1,Y_2]$ given by $$\Psi=\Big{[}aX,\mu Y_1+\tilde c(X),a^{-d}{\mu}^rY_2+(aX)^{-d}\Big{(}{\underset{i=1}{\overset{r}\prod} \big{(}\mu Y_1+\tilde c(X)-\sigma_i(aX)\big{)}-{\mu}^r{\underset{i=1}{\overset{r}\prod} \big{(}Y_1-\sigma_i(X)}\big{)}}\Big{)}\Big{]},$$ where $\tilde c(X)=c(X)+X^db(X)$. It follows that the composition $\Phi_2\circ\Phi_1$ of two automorphisms $\Phi_1$ and $\Phi_2$ of $B_{d,P}$ with data $\mathcal{A}_{\Phi_1}=\big{(}\alpha_1,\mu_1,a_1, b_1(x)\big{)}$ and $\mathcal{A}_{\Phi_2}=\big{(}\alpha_2,\mu_2,a_2, b_2(x)\big{)}$  respectively is the automorphism of $B_{d,P}$ with datum $$\mathcal{A}=\big{(}\alpha_2\alpha_1,\mu_2\mu_1,a_2a_1,{a_2}^{-d}\mu_2b_1(x)+b_2(a_1x)\big{)}. \quad\quad\quad\quad~~~~~~~~(*)$$

\noindent 
 In Theorem 3.11 of \cite{DP}, it has been shown that the following statements hold.
 \begin{enumerate}

 \item[\rm(i)] For each $b(x) \in k[x]$, the datum $\big{(}id,1,1,b(x)\big{)}$ induces an automorphism of $B_{d,P}$ given by $U_{b(x)}=\Big{[}x, y_1+x^db(x), y_2+x^{-d}\Big{(}P\big{(}x,y_1+x^db(x)\big{)}-P(x,y_1)\Big{)}\Big{]}$.
 
 \item[\rm(ii)] If $a^q\neq 1$ for every $q\in{\{1,\dots,d-1\}}$, then there exist $\tau(x) \in k[x]$, $Q(y_1) \in k[y_1]$ such that $P\big{(}x,y_1+\tau(x)\big{)}=Q(y_1)$ and the datum $(id,1,a,0)$ induces an automorphism of $B_{d,P}$ given by $H_a=\Big{[}ax, y_1+\tau(ax)-\tau(x), a^{-d}y_2\Big{]}$.
 
 \item[\rm(iii)] If $a^{q_0}=1$ for some minimal $q_0\in{\{2,\dots,d-1\}}$, then there exist $\tau(x) \in k[x]$, $Q(x,y_1) \in k[x,y_1]$ such that $P\big{(}x,y_1+\tau(x)\big{)}=Q(x^{q_0},y_1)$ and the datum $(id,1,a,0)$ induces an automorphism of $B_{d,P}$ given by $H_{a,q_0}=\Big{[}ax, y_1+\tau(ax)-\tau(x), a^{-d}y_2\Big{]}$.
 
 \item[\rm(iv)]If $\alpha=id$, then $\mu=1$. If $\alpha\neq{id}$, then there exists integer $s\geqslant 2$ such that $\mu^s=1$ and $\mu^{s^{\prime}}\neq 1$ for $s^{\prime}\in{\{1,\dots,s-1\}}$. Moreover, if $\alpha\neq id$, then there exist $\tau(x) \in k[x]$, $Q(x,y_1) \in k[x,y_1]$ such that $P\big{(}x,y_1+\tau(x)\big{)}=y_1^iQ(x,y_1^s)$ {\rm(}where $i=0,1)$ and the datum $(\alpha, \mu ,1, 0)$ induces an automorphism of $B_{d,P}$ given by $S_{\alpha,\mu}=\Big{[}x, \mu y_1+(1-\mu)\tau(x), {\mu}^ iy_2\Big{]}$.
 \end{enumerate}

\noindent
Let ${{\U}},{\bH},{\bS}$ be the subgroups of ${\rm Aut}(B_{d,P})$ consisting of the automorphisms corresponding to the data of the type $\big{(}id,1,1,b(x)\big{)},(id,1,a,0)$ and $(\alpha, \mu, 1,0)$ respectively. Since ${\rm ML}(B_{d,P})=k[x]$ (\cite[Theorem 2.6]{D}), it follows from the statements above, that $B_{d,P}$ is an almost rigid domain with the canonical {\it lnd} $D_{d,P}$, given by
$$D_{d,P}:=x^d\dfrac{\partial}{\partial y_1}+\dfrac{\partial P}{\partial y_1} \dfrac{\partial}{\partial y_2}.$$ 

\begin{rem}\label{dansaut}
{\em $\U$ is the big unipotent group $\U(D_{d,P})$ and the automorphisms ${H}_a$, $H_{a,q_0}$ correspond to a non-trivial action of an one-dimensional algebraic torus $k^*$ on $B_{d,P}.$
}
\end{rem}  

With the notation as above one obtains the following structure of ${\rm Aut}(B_{d,P})$.
\begin{lem}\label{Aut_DS}
${\rm Aut}(B_{d,P})\cong{({\bS}\times{{\bH}})}\ltimes{\U}$.
\end{lem}
\begin{proof}
 Let $\Phi$ be an automorphism of $B_{d,P}$ corresponding to the datum $\big{(}\alpha,\mu,a,b(x)\big{)}$. If $\alpha=id$, then by (iv) we have  $\mu=1$. By ($*$), we observe that $$\big{(}id,1,a,b(x)\big{)}=\big{(}id,1,a,0\big{)}\circ \big{(}id,1,1,b(x/a)\big{)}.$$ Hence we can write $\Phi$ as
 $$\Phi=H\circ U_{b(x/a)},$$ where $H=
 \begin{cases}
  Id_{B_{d,P}},& \text{if }a=1.\\
  H_a,& \text{if }a^q\neq 1$ for any $q\in{\{1,\dots,d-1\}}.\\
  H_{a,q_0},&\text{if }a^q\neq 1$ for every $q\in{\{1,\dots,d-1\}}.
 \end{cases}
$

\noindent\smallskip
If $\alpha\neq id$, then $\mu\neq 1$. It follows from ($*$) that $$\big{(}\alpha,\mu,a,b(x)\big{)}=\Big{(}\big{(}\alpha,\mu,1,0\big{)}\circ\big{(}id,1,a,0\big{)}\Big{)}\circ \big{(}id,1,1,b(x/a)\big{)}.$$ Hence we can write $\Phi$ as
 $$\Phi=(S_{\alpha,\mu}\circ H)\circ U_{b(x/a)},$$ where $H$ is as above. So, ${\rm Aut}(B_{d,P})=\bS\bH\U$.

 \noindent\smallskip
 Since the elements of $\bS$ and $\bH$ commute and $\bS\cap\bH=\{Id_{B_{d,P}}\}$, we have $\bS\bH\cong{\bS\times\bH}$. Moreover, $\U$ is a normal subgroup of ${\rm Aut}(B_{d,P})$ and $\U\cap\bS\bH=\{Id_{B_{d,P}}\}$. Hence the result follows.
\end{proof}

\indent
We now prove a proposition which give necessary and sufficient conditions for an element of ${\rm Aut}(B_{d,P})$ to be in the isotropy subgroup of any {\it lnd} on $B_{d,P}$.
\begin{lem}\label{y_1DS}
Let $\delta\in{\rm LND}(B_{d,P})$ and $\alpha\in{\rm Aut}(B_{d,P})$. Then  the following are equivalent.
\begin{enumerate}
\item[\rm(i)]$\alpha\in {\rm Aut}(\delta)$.
\item[\rm(ii)]$\delta\alpha(y_1)=\alpha \delta(y_1)$.
\item[\rm(iii)]$\delta\alpha(y_2)=\alpha \delta(y_2)$.
\end{enumerate}
\end{lem}
\begin{proof}
Since  $x^dy_2=P(x,y_1)$ and by Theorem \ref{Aut_DS}, $\alpha(x)=ax$ for some $a\in{k^*}$, we have
\begin{enumerate}
\item[\rm(1)] $a^{d}x^{d}\alpha(y_2)=P\big{(}\alpha(x),\alpha(y_1)\big{)}$,

\item[\rm(2)] $x^{d}\delta(y_2)=\dfrac{\partial P}{\partial {y_1}}\delta(y_1)$.

\end{enumerate}

\noindent
Since $\delta\alpha(x)=0=\alpha \delta(x)$, it is enough to show that $\rm(ii)\Leftrightarrow \rm(iii)$.

\medskip
\noindent
Now, $$\alpha \delta(y_2)=\delta\alpha(y_2)\Leftrightarrow{a^{d}x^{d}\alpha \delta(y_2)=a^{d}x^{d}\delta\alpha(y_2)}\Leftrightarrow{\alpha\rbo x^{d}\delta(y_2)\rbc=\delta\rbo a^{d}x^{d}\alpha(y_2)\rbc}$$ Hence by (1) and (2), $$\alpha \delta(y_2)=\delta\alpha(y_2)\Leftrightarrow{\alpha\rbo\dfrac{\partial P}{\partial {y_1}}\delta(y_1)\rbc=\delta\Big{(}P\big{(}\alpha(x),\alpha(y_1)\big{)}\Big{)}}.$$ Therefore by Lemma \ref{chainrule}, $$\alpha \delta(y_2)=\delta\alpha(y_2) \Leftrightarrow{\alpha(\dfrac{\partial P}{\partial {y_1}})\alpha \delta(y_1)=\alpha(\dfrac{\partial P}{\partial {y_1}})\delta\alpha(y_1)}.$$

\smallskip
\noindent
Since $\dfrac{\partial P}{\partial {y_1}}\neq{0}$, $\alpha \delta(y_2)=\delta\alpha(y_2) \Leftrightarrow{\alpha \delta(y_1)=\delta\alpha(y_1)}$. 

\end{proof}

The following theorem is the main result of this section.
\begin{thm}\label{thmdp}
 Let $\delta=f(x)D_{d,P}$, where $f(x)={\underset{i=0}{\overset{l}\sum}}{a_ix^{n_i}}\in{k[x]}~(n_i,l \in \bN \cup \{0\}, a_i \in k^*$ for each $i)$. Suppose $\Ga:={\rm Aut}(\delta)$ and $n:={\rm GCD}(d+n_0, \dots, d+n_l)$. 
 Then the following statements hold.
\begin{enumerate}
\item[\rm(i)] The unipotent group ${\U}\subseteq{\Ga}$ and hence ${\U}_{\delta}=\U$.

\item[\rm(ii)] If $a \in k^*$ with $a^q\neq 1$ for any $q\in{\{1,\dots,d-1\}}$, then ${H}_a\in {\Ga}$ if and only if  $a^n=1$ and $n\geqslant d$.

\item[\rm(iii)] If $a \in k^*$ with $a^{q_0}=1$ for some minimal $q_0\in{\{2,\dots,d-1\}}$, then ${H}_{a,q_0}\in {\Ga}$ if and only if $q_0\mid{n}$.

\item [\rm (iv)] If $k$ is algebraically closed, then $\bH_{\delta} \cong \bZ_n$. Moreover, if $n\geqslant d$, then the generators of $\bH_{\delta}$ are of the form $H_a$, where $a$ is a primitive $n$-th root of unity. Otherwise, they are of the form $H_{a,n}$.

\item[\rm(v)] The subgroup ${\bS}_{\delta}=\{{ Id}_{B_{d,P}}\}$.
\item[\rm(vi)] The isotropy subgroup $\Ga\cong{{\U}}\rtimes{({\bH}\times{{\bS}})_{\delta}}.$

\end{enumerate} 
\end{thm}

\begin{proof}

\begin{enumerate}
\item[\rm(i)] Let $b(x) \in k[x]$. Then, $$\delta{\U}_{b(x)}(y_1)=\delta\big{(}y_1+x^db(x)\big{)}=f(x)x^d$$
and
$${\U}_{b(x)}\delta(y_1)=\alpha_b\big{(}f(x)x^d\big{)}=f(x)x^d.$$
Hence, by Lemma \ref{y_1DS}, ${\U}_{b(x)}\in {\rm Aut}(\delta)$. 
\item[\rm(ii)] Let $a \in k^*$ be such that $a^q\neq 1$ for any $q\in{\{1,\dots,d-1\}}$. Then,
$$\delta {H}_a(y_1)=\delta\rbo y_1+\tau(ax)-\tau(x)\rbc=f(x)x^d$$ and
$${H}_a\delta(y_1)=H_a\rbo f(x)x^d\rbc=a^dx^df(ax).$$ By Lemma \ref{root of unity}, ${H}_a \in {\rm Aut}(\delta)$ if and only if $a^n=1$  and $n\geqslant d$.
\item[\rm(iii)] Let $a \in k^*$ be such that $a^{q_0}=1$ for some minimal $q_0\in{\{2,\dots,d-1\}}$. Then, $${H}_{a,q_0} \in {\rm Aut}(\delta) \Leftrightarrow a^{m}=1, \text{where}~m={\rm GCD}(n,q_0)~ \text{and}~m\geqslant {q_0} \Leftrightarrow q_0\mid n.$$
\item[\rm(iv)] The proof follows from (ii) and (iii).
\item[\rm(v)] Let $\alpha\in {\mathcal S}_r$ and $\mu\in{k^*}$. If $\alpha=id,$ then $\mu=1$ and hence the datum $(\alpha,\mu,1,0)$ corresponds to $Id_{B_{d,P}}.$ If $\alpha\neq id$, then $\mu\neq 1$ and we observe that $$\delta {S}_{\alpha,\mu}(y_1)=\delta\rbo\mu y_1+(1-\mu)\tau(x)\rbc=\mu f(x)x^d$$ and

$${S}_{\alpha,\mu} \delta(y_1)=S_{\mu}\rbo f(x)x^d \rbc=f(x)x^d.$$
Hence,  $\bS_{\delta}=\{Id_{B_{d,P}}\}$.

\item[\rm(vi)]The proof follows from Lemma \ref{Aut_DS}.

\end{enumerate}
\end{proof}

\begin{rem}\label{dansrem} {\em When $k$ is algebraically closed, from Theorem \ref{thmdp}, we can immediately make the following observations.

\begin{enumerate}
\item [\rm (i)] It may be possible that $({\bH}\times {\bS})_{\delta} \supsetneqq{\bH}_{\delta}\times{\bS}_{\delta}$. Let $H=[ax, y_1+\tau_1(ax)-\tau_1(x), a^{-d}y_2] \in {\bH}$ and $S=[x, \mu y_1+(1-\mu)\tau_2(x), {\mu}^ iy_2] \in {\bS}$.
Then ${H}{S}=[ax, \mu y_1+ (1-\mu)\tau_2(x)+\tau_1(ax)-\tau_1(x), a^{-d}\mu^iy_2]$. 
Now, $${H}{S}\delta(y_1)={H}{S}\big{(}f(x)x^d\big{)}=a^dx^df(ax)$$ and $$\delta{H}{S}(y_1)=\mu x^df(x).$$ Hence, ${H}{S}\in ({\bH}\times{\bS})_{\delta}$ if and only if $\mu=a^{d+n_i}$ for all $i$, $0\leqslant i\leqslant l$. So, if we choose $f(x)$ to be in $k^*$, $\mu\in{k^*}$ be such that $\mu\neq 1$ and $a$ to be a $d$-th root of $\mu$, then neither ${H}\notin{{\bH}}_{\delta}$ nor ${S}\notin{{\bS}}_{\delta}$.
\item [\rm (ii)] If $\delta=D_{d,P}$, then $\U(\delta) \subsetneqq {\rm Aut}(\delta)$. Indeed, if $\omega \in k$ be a primitive $d$-th root of unity, then ${H}_{\omega} ~(\neq Id_{B_{d,P}}) \in {\rm Aut}(\delta)$. However, it may happen that ${\rm Aut}(\delta)=\U$, when $\delta$ is a replica of $D_{d,P}$. For example, if $\delta=(x+x^2)D_{d,P}$, then $({{\bH} \times {\bS})}_{\delta}=\{Id_{B_{d,P}}\}$ as $n=1$.
\end{enumerate}
}
\end{rem}

\section{Danielewski varieties with constant coefficients}\label{iso_Dvc}

\indent
For an affine algebraic variety $\X$ over $\bK$, we will denote the coordinate ring or equivalently, the algebra of regular functions of $\X$ by $\bK[\X]$ and the group of automorphisms of $\bK[\X]$ by ${\rm Aut}(\X)$. For convenience, we also denote ${\rm LND}(\bK[\X])$ by ${\rm LND}(\X)$.
The following generalization of {\it Danielweski surfaces} were introduced by Dubouloz (\cite{D}) as counter-examples to the generalized Zariski cancellation problem in arbitrary dimension. Throughout this section, $d,m,k_2,\dots,k_m$ are natural numbers greater or equal to $2$.

\begin{defn}\label{danvar} An affine algebraic variety $\V \subseteq {\bK}^{m+1}$ is called a Danielewski variety if 
$$\bK[\V]=\dfrac{\bK[Y_1,Y_2, \dots ,Y_m,Z]}{(Y_1Y_2^{k_2} \dots Y_m^{k_m}-P(Y_2, \dots, Y_m,Z)},$$
where 
$P(Y_2, \dots ,Y_m,Z)=Z^d+s_{d-1}(Y_2,\dots ,Y_m)Z^{d-1}+ \dots + s_0(Y_2, \dots ,Y_m).$
\end{defn}

\noindent
If $s_i\in \bK$ for all $i\in\{0,\dots,d-1\}$, then $\V$ will be called a {\it Danielewski variety with constant coefficients}, denoted by $\V_{con}$.

\smallskip
\noindent
Let $y_1,y_2, \dots ,y_m,z$ denote the images of $Y_1,Y_2, \dots , Y_m,Z$ respectively in $\bK[\V_{\it con}]$. In Section 6 of \cite{Gai}, the locally nilpotent derivations on $\V_{con}$ have been studied and the structure of ${\rm Aut}(\V_{con})$ have been described in detail. The variety $\V_{con}$ is irreducible (\cite[Lemma 4.1]{Gai}) with no non-constant invertible functions (\cite[Lemma 4.2]{Gai}) and $\bK[\V_{\it con}]$ is an almost rigid domain (\cite [Lemma 6.2]{Gai}) with the {\it canonical lnd} $D_{con}$, given by 
$$ D_{con}(y_1)=P^{\prime}(z),~~D_{con}(z)=y_2^{k_2}\dots y_m^{k_m}~~\text{and}~~D_{con}(y_j)=0,~\forall j \geqslant 2.$$

\smallskip
\noindent
Consider the natural diagonal action of the $m$-dimensional algebraic torus ${({\bK}^*)}^m$  on $\bK[Y_1,Y_2, \dots ,Y_m]$. The stabilizer of the monomial $Y_1Y_2^{k_2}\dots Y_m^{k_m}$ is isomorphic to the $(m-1)$-dimensional torus, denoted by $\bT$. If we consider the trivial action of $\bT$ on $\bK[Z]$, then there is an effective action of $\bT$ on $\V_{con}$. $\bT$ is called the {\it proper torus} of $\V_{con}$.

\smallskip
\noindent
There is a natural action of the symmetric group $\mathcal{S}_m$ on $\bK[Y_1,\dots ,Y_m]$. The stabilizer $\bS$ of the monomial $Y_1Y_2^{k_2}\dots Y_m^{k_m}$ is isomorphic to the group
$$\mathcal S_{m_1} \times \dots \times \mathcal S_{m_n},$$ where $m=m_1+\dots +m_n$ and for each $i\in \{1,\dots,n\}$, $\mathcal S_{m_i}$ permutes the $m_i$ many $Y_j$'s with same $k_j$. If we consider the trivial action of $\bS$ on $\bK[Z]$, then there is an effective action of $\bS$ on $\V_{con}$. This group $\bS$ is called the {\it symmetric group} of $\V_{con}$. 

\smallskip
\noindent
If $P(Z)=Z^d$, there is also an effective action of an one-dimensional torus ${\bK}^*$ acting by $$t\cdot (y_1,y_2, \dots,y_m,z)=(t^dy_1,y_2 \dots ,y_m,tz),~~\text{for all}~t \in {\bK}^*.$$ If $P(Z) \neq Z^d$ and $v$ is the maximal integer such that there exists a polynomial $Q(Z) \in \bK[Z]$ and a nonnegative integer $u$ such that $P(Z)=Z^uQ(Z^v)$, then there is an action of ${\bZ}_v$ (considered as a subgroup of $k^*$) on $\V_{con}$, given by $$t \cdot (y_1,y_2,\dots,y_m,z)=(t^uy_1,y_2\dots ,y_m,  tz),~~\text{where}~t^v=1.$$  In each of these cases, the groups ${\bK}^*$ and $\bZ_v$ are called the {\it additional quasitorus} of $\V_{con}$ and denoted by $\bD$. One can refer to Section 4 of \cite{Gai} for further details. 

\smallskip
\noindent
Let us denote the subgroups of ${\rm Aut}(\V_{con})$ generated by the automorphisms induced by the proper torus, symmetric group and additional quasitorus by $\bT$, $\bS$ and $\bD$ respectively.

\smallskip
\noindent
\begin{prop}\label{ptorus1}
Let $\delta =hD_{con}$, where $h \in \bK[y_2,\dots ,y_m]$ and $\theta \in {\bT}$. Then the following statements are equivalent.
\begin{enumerate}
\item [\rm (a)] $\delta \theta (y_1)=\theta \delta (y_1)$.
\item [\rm (b)] $h\theta (y_1)=y_1 \theta (h)$.
\item [\rm (c)] $\delta \theta (z)=\theta \delta (z)$.
\item [\rm (d)] $\theta \in {\rm Aut}(\delta)$.
\end{enumerate}
Moreover, if $\theta(h)=h$ $($in particular, when $h \in {\bK}^*)$, the statement ${\rm(b)}$ is equivalent to 

\begin{enumerate}
\item[$\rm(b)^{\prime}$] $ \theta(y_1)=y_1$.
\end{enumerate}

\end{prop}

\begin{proof}

\smallskip\noindent
{\bf \underline{(a)$\Leftrightarrow$(b)}}

\smallskip
\noindent
Let $\theta (y_1)={\lambda}_1y_1$, for some $\lambda_1 \in {\bK}^*$. Since $\bT$ acts trivially on $\bK[z]$, we have $\theta(z)=z$. Then,
$$\delta \theta(y_1)=\lambda_1hP^{\prime}(z)~~\text{and}~~\theta \delta(y_1)=\theta(h)P^{\prime}(z).$$ So,
$$\delta \theta (y_1)=\theta \delta (y_1) \Leftrightarrow h\lambda_1= \theta (h).$$
Since $\V_{con}$ is irreducible, $h\lambda_1= \theta (h)\Leftrightarrow h\theta (y_1)=y_1 \theta (h)$.

\smallskip
\noindent
{\bf \underline{(b)$\Rightarrow$(c)}}
We have \begin{equation*}
\begin{split}
y_1\delta \theta(z) & =hy_1y_2^{k_2}\dots y_m^{k_m} \\
&=h\theta(y_1y_2^{k_2}\dots y_m^{k_m})~~~~~~~~~~\text{(since}~\theta~\text{stabilizes}~y_1y_2^{k_2}\dots y_m^{k_m})\\
& = \big{(}h \theta (y_1)\big{)}\theta (y_2^{k_2}\dots y_m^{k_m})\\
& = \big{(}\theta (h) y_1\big{)}\theta (y_2^{k_2}\dots y_m^{k_m})~~~~~~~~~\text{\big{(}by (b)\big{)}}\\
& = y_1 \theta (hy_2^{k_2}\dots y_m^{k_m})=y_1\theta \delta(z).
\end{split}
\end{equation*}
Since $\V_{con}$ is irreducible, (b)$\Rightarrow $(c). 

\smallskip
\noindent
{\bf \underline{(c)$\Rightarrow$(b)}}

\smallskip
\noindent
Suppose $\delta \theta (z)=\theta \delta (z)$.
$$\text{Now,}~~y_1\theta(h)\theta (y_2^{k_2}, \dots ,y_m^{k_m})=y_1\theta\big{(}\delta (z)\big{)}=y_1\delta \theta (z)=y_1hy_2^{k_2}, \dots ,y_m^{k_m}=h\theta(y_1y_2^{k_2}, \dots ,y_m^{k_m}).$$ Since $\V_{con}$ is irreducible, $y_1\theta(h)=h\theta(y_1)$.

\smallskip
\noindent
{\bf \underline{(a)$\Leftrightarrow$(d)}}

\smallskip
\noindent
Clearly, (d)$\Rightarrow$(a). Since (a), (c) are equivalent and for any $j \in \{2,\dots ,m\}$, 
$ \delta \theta (y_j) =0=\theta\delta (y_j)$, hence (a)$\Rightarrow$(d). 

\smallskip
\noindent
If $\theta(h)=h$, then the equivalence of (b) and {(b)}$^{\prime}$ is obvious.
\end{proof}

\smallskip
The next corollary gives the structure of ${\bT}_{\delta}$ when $\delta$ is a scalar multiple of $D_{con}$.
\begin{cor}\label{ptorus2}
Let $\delta=hD_{con}$, for some $h \in {\bK}^*$. Then ${\bT}_\delta\cong{(\bK^*)}^{m-2}\times\bZ_s$, where $s={\rm GCD}(k_2,\dots,k_m)$.
\end{cor}

\begin{proof} Recall that $\bT \cong {(\bK^*)}^{m-1}$. Let $\theta \in {\bT}_{\delta}$. Then, for each $j \in \{1,2,\dots ,m \}$, there exists $\lambda_j \in {\bK}^*$ such that $\theta (y_j)=\lambda_jy_j$. Since $\theta$ stabilizes $y_1y_2^{k_2}\dots y_m^{k_m}$, we have $\lambda_1=\frac{1}{\lambda_2^{k_2}\dots \lambda_m^{k_m}}.$ Since $\theta(h)=h$, by Proposition \ref{ptorus1}, we have $\lambda_1=1$. So, ${\bT}_{\delta}$ is isomorphic to the subgroup $$\{(\lambda_2, \dots ,\lambda_m) \in {({\bK}^*)}^{m-1}~|~\prod_{j=2}^{m}\lambda_j^{k_j}=1\}.$$ Hence, the result follows.
\end{proof}

\begin{rem} {\em If $\delta =hD_{con}$, for some $h \in {\bK}^*$ (in particular, when $\delta=D$), then, by Corollary \ref{ptorus2}, ${\bT}_{\delta} \neq \{Id\}$. However, the next example shows that when $h \in \bK[y_2,\dots,y_m]\setminus {\bK}^*$, it may be possible that ${\bT}_{\delta}=\{Id\}$.
}
\end{rem}
\begin{ex}\label{ptorusex1}
{\em Let $h={y_2}^2+y_2+\dots+y_m$. Suppose there exists $\theta \in{\bT_\delta} $ and $\theta(y_j)=\lambda_jy_j$, for $\lambda_j\in{\bK^*}$, $j=1,2,\dots,m$. Hence, by Proposition \ref{ptorus1}, $\theta(h)=\lambda_1h$, which is true if and only if $\lambda_j=1$ for all $j \geqslant 1$. Hence, $\theta=Id$.
}
\end{ex}

\smallskip
\indent
The next lemma gives a necessary and sufficient condition for an element of $\bS$ to be in ${\rm Aut}(\delta)$, for any $\delta\in{\rm LND}(\V_{con})$.
\begin{lem}\label{symmetric}
Let $\delta =hD_{con}$, for some $h \in \bK[y_2, \dots ,y_m]$ and $\sigma \in {\bS}$. Then the following statements are equivalent.
\begin{enumerate}
\item [\rm (a)] $\delta \sigma (y_1)=\sigma \delta (y_1)$.
\item [\rm (b)] $\sigma (h)=h$.
\item [\rm (c)] $\delta \sigma (z)=\sigma \delta (z)$.
\item [\rm (d)] $\sigma \in {\rm Aut}(\delta)$.
\end{enumerate}
In particular, if $h \in {\bK}^*$, then $\bS~(={\bS}_{\delta})$ is a subgroup of ${\rm Aut}(\delta)$.
\end{lem}

\begin{proof} 
 We first make the following observations.

\begin{enumerate}
	\item [\rm(i)] Since each $k_j \geqslant 2$, for any $\sigma \in \bS$, we have $\sigma (y_1)=y_1$.
	\item [\rm(ii)]  Since $\sigma$ stabilizes both $y_1$ and $y_1y_2^{k_2} \dots y_m^{k_m}$ and $\V_{con}$ is irreducible, we have $\sigma(y_2^{k_2} \dots y_m^{k_m})=y_2^{k_2} \dots y_m^{k_m}$.
\end{enumerate}
{\bf \underline{(a)$\Leftrightarrow$(b)}}

\smallskip\noindent
$\delta \sigma (y_1)=\delta (y_1)=hP^{\prime}(z),$ and
$\sigma \delta (y_1)=\sigma (hP^{\prime}(z))=\sigma (h) P^{\prime}(z).$
Hence, (a)$\Leftrightarrow(b)$.

\smallskip
\noindent
{\bf \underline{(b)$\Leftrightarrow$(c)}}

\smallskip\noindent
Suppose (b) holds. Then
$$
y_1\delta \sigma (z)=y_1\delta (z)=hy_1y_2^{k_2} \dots y_m^{k_m} 
=y_1 \sigma (h) \sigma (y_2^{k_2} \dots y_m^{k_m})=y_1 \sigma \delta (z).$$

Hemce, $\delta \sigma (z)=\sigma \delta (z)$. 

\noindent
Conversely, suppose (c) holds. Then,
$$hy_2^{k_2} \dots y_m^{k_m}=\delta \sigma(z)=\sigma \delta (z)=\sigma (h y_2^{k_2} \dots y_m^{k_m})=\sigma (h) y_2^{k_2} \dots y_m^{k_m}.$$
Hence, $\sigma(h)=h$.

\smallskip
\noindent
{\bf \underline{(a)$\Leftrightarrow$(d)}}

\smallskip\noindent
Clearly, (d)$\Rightarrow$(a). Conversely, suppose (a) holds. Since $\sigma$ permutes $\{y_2,\dots,y_m\}$, 
$\delta \sigma (y_j)=0=\sigma \delta (y_j)$, for any $j \geqslant 2$ and (a)$\Rightarrow$(c), we have $\sigma \in {\rm Aut}(\delta)$.

\smallskip
\noindent
Moreover, if $h \in {\bK}^*$, then $\sigma (h)=h$, for all $\sigma \in \bS$ and $h \in \bK[y_2, \dots y_m]$. Hence the result follows.
\end{proof}

\smallskip
\indent
The next result shows that the no non-trivial automorphism in $\bD$ can be in ${\rm Aut}(\delta)$, for any $\delta \in {\rm LND}(\V_{con})$. 
\begin{lem}\label{qtorus}
Let $\delta =hD_{con}$, for some $h \in \bK[y_2, \dots ,y_m]$. Let $\varphi \in \bD$. Then $$\varphi \in {\rm Aut}(\delta) \Leftrightarrow \varphi=Id.$$
\end{lem}

\begin{proof}
Since $\varphi\in{\bD}$, there exists $t\in {\bK}^*$ such that 
$$\varphi(y_1,y_2,\dots,y_m,z)=
\begin{cases}
(t^dy_1,y_2,\dots,y_m,tz), & \text{if } P(Z)=Z^d\\
(t^uy_1,y_2,\dots,y_m,tz), & \text{if } P(Z)=Z^uQ(Z^v).
\end{cases}
$$ If $\varphi \in {\rm Aut}(\delta)$, then

$$
thy_2^{k_2}\dots y_m^{k_m}=\delta \varphi (z)= \varphi \delta (z)= \varphi (hy_2^{k_2} \dots y_m^{k_m})= hy_2^{k_2} \dots y_m^{k_m}.$$

\noindent
So, $t=1$ and hence $\varphi=Id$. The converse is obvious.
\end{proof}

\smallskip
\indent
Let $\mathbb{U}(D_{con})$ be the big unipotent subgroup of ${\rm Aut}(\V_{con})$ corresponding to the {\it canonical lnd} $D_{con}$ on $\V_{con}$. Clearly, $\mathbb{U}(D_{con})$ is a subgroup of ${\rm Aut}(\delta)$ for every $\delta \in {\rm LND}(\V_{con})$. The structure of ${\rm Aut}(\V_{con})$ has been described in \cite{Gai} (c.f. Proposition 6.3 and Theorem 6.4). We state the precise result below.
\begin{prop}\label{danst1}
 ${\rm Aut}(\V_{con})\cong\bS \ltimes \big{(}(\bT \times \bD) \ltimes \mathbb{U}(D_{con}) \big{)}.$
\end{prop}

We now state the main result of this section.
\begin{thm}\label{danst2}
The isotropy subgroup ${\rm Aut}(D_{con})$ can be described as follows.
\begin{enumerate}

\item[\rm(i)]If $P(Z)=Z^d$, then $${\rm Aut}(D_{con})\cong{{\bS} \ltimes \big{(}{(\bK^*)}^{m-1} \ltimes \mathbb{U}(D_{con})\big{)}}.$$ 
 
\item[\rm(ii)]If $P(Z)\neq Z^d$ and $v$ is the maximal integer such that $P(Z)=Z^uQ(Z^v)$, then $${\rm Aut}(D_{con})\cong{\bS} \ltimes \Big{(}\big{(}{(\bK^{*})}^{m-2}\times {\bZ}_{sv}\big{)}\ltimes \mathbb{U}(D_{con})\Big{)},$$ where $s={\rm GCD}(k_2,\dots,k_m)$.
 \end{enumerate}
\end{thm}

\begin{proof} 
Let $\sigma\in{\rm Aut}(D_{con})$. Then, by Proposition \ref{danst1}, $\sigma=\sigma_1\sigma_2\sigma_3\sigma_4$, for some $\sigma_1\in{\bS}$, $\sigma_2\in{\bT}$, $\sigma_3\in{\bD}$ and $\sigma_4\in{\U(D_{con})}$. Since $\U(D_{con})\subseteq{\rm Aut}(\delta)$ and ${\bS}_{D_{con}}=\bS$ (by Lemma \ref{symmetric}), we have $\sigma_2\sigma_3\in{({{\bT}\times{\bD}})_{D_{con}}}$. 

\smallskip\noindent
(i) Let $\sigma_2(y_1)=\lambda y_1$ and $\sigma_3(z)=tz$, for some $t,\lambda\in{\bK^{*}}$. Then
$$\delta\sigma_2\sigma_3(y_1)=\delta\sigma_2(t^dy_1)=\delta(t^d\lambda y_1)=t^d\lambda P^{\prime}(z) $$ and 
$$\sigma_2\sigma_3\delta(y_1)=\sigma_2\sigma_3\big{(}P^{\prime}(z)\big{)}=\sigma_2\big{(}t^{d-1}P^{\prime}(z)\big{)}=t^{d-1}P^{\prime}(z).$$
So, $\sigma_2\sigma_3\in{({{\bT}\times{\bD}})_{D_{con}}}$ if and only if $t\lambda=1$. Hence, $({{\bT}\times{\bD})}_{D_{con}}\cong{{(\bK^{*})}^{m-1}}$.

\smallskip
\noindent
(ii) Let $\sigma_2(y_i)=\lambda_i y_i$, for $i=1,\dots,m$ and  $\sigma_3(z)=tz$, where $\lambda_i,t\in{\bK^{*}}$ with\\ $\lambda_1=\dfrac{1}{\prod_{i=2}^{m} \lambda_i^{k_i}}$ and $t^v=1$. Then $$\delta\sigma_2\sigma_3(y_1)=\delta\sigma_2(t^uy_1)=\delta(t^u\lambda_1 y_1)=t^u\lambda_1 P^{\prime}(z) $$ and
$$\sigma_2\sigma_3\delta(y_1)=\sigma_2\big{(}t^{u-1}P^{\prime}(z)\big{)}=t^{u-1}P^{\prime}(z).$$
Hence, $$\sigma_2\sigma_3\in{({{\bT}\times{\bD}})_{D_{con}}}\Leftrightarrow t\lambda_1=1\Leftrightarrow t=\prod_{i=2}^{m} \lambda_i^{k_i}.$$  
So, 
\begin{equation*}
\begin{split}
({{\bT}\times{\bD})}_{D_{con}} &\cong\{(\dfrac{1}{\prod_{i=2}^{m} \lambda_i^{k_i}},\lambda_2,\dots,\lambda_m,\prod_{i=2}^{m} \lambda_i^{k_i})\in{(\bK^*)}^{m+1}:\prod_{i=2}^{m} \lambda_i^{vk_i}=1\}\\
&\cong{\{(\lambda_2,\dots,\lambda_m)\in{(\bK^*})^{m-1}:\prod_{i=2}^{m} \lambda_i^{vk_i}=1}\}\\
&\cong{{{(\bK^{*})}^{m-2}\times{\bZ}_{sv}}}
\end{split}
\end{equation*}
\end{proof}

%
%
%
%
%
%
%
%
%

\begin{rem}\label{dancrem} {\em
From Theorem \ref{danst2}, we can  make the following observations. 
\begin{enumerate}

\item[\rm(i)] If $P(Z)=Z^d$, then there exist $\sigma_1\in\bT$ and $\sigma_2\in\bD$ such that $\sigma_1\sigma_2\in{{(\bT\times\bD)}_{D_{con}}}$ but neither $\sigma_1\in\bT_{D_{con}}$ nor $\sigma_2\in\bD_{D_{con}}~(=\{Id\},~\text{by Lemma \ref{qtorus}})$. For example, let $$\sigma_1(y_1,\dots,y_m,z)=(\dfrac{1}{\lambda^{k_2}}y_1,\lambda y_2,y_3,\dots,y_m,z)~~\text{and}$$ $$\sigma_2(y_1,y_2,\dots,y_m,z)=(\lambda^{k_2d}y_1,y_2,\dots,y_m,\lambda^{k_2} z),$$ where $\lambda\in\bK^*$ be such that $\lambda^{k_2}\neq 1$.
\item[\rm(ii)] If $P(Z)\neq Z^d$ and $v\geqslant 2$ be the maximal integer such that $P(Z)=Z^uQ(Z^v)$ then there exist $\sigma_1\in\bT$ and $\sigma_2\in\bD$ such that $\sigma_1\sigma_2\in{{(\bT\times\bD)}_{D_{con}}}$ but neither $\sigma_1\in\bT_{D_{con}}$ nor $\sigma_2\in\bD_{D_{con}}~(=\{Id\})$. For example, let $$\sigma_1(y_1,\dots,y_m,z)=(\dfrac{1}{\lambda^{k_2}}y_1,\lambda y_2,y_3,\dots,y_m,z) ~~\text{and}$$ $$\sigma_2(y_1,y_2,\dots,y_m,z)=(\lambda^{k_2u}y_1,y_2,\dots,y_m,\lambda^{k_2} z),$$ where $\lambda\in\bK^*$ be such that $\lambda^{k_2v}=1$ but $\lambda^{k_2}\neq 1$.
\item [\rm (iii)] If $P(Z)=Z^d$, then ${\rm Aut}(\V_{con})$ contains a subgroup isomorphic to the $(m-1)$-dimensional  algebraic torus and hence properly contains $\U(D_{con})$. If $P(Z)\neq Z^d$ and $v$ is the maximal integer such that $P(Z)=Z^uQ(Z^v)$, then ${\rm Aut}(\V_{con})$ contains a subgroup isomorphic to an $(m-2)$-dimensional algebraic torus whenever $m\geqslant 3$ and for $m=2$, ${\rm Aut}(\V_{con})$ contains a subgroup isomorphic to $\bZ_{k_2v}$. Hence, ${\rm Aut}(\V_{con})$ properly contains $\U(D_{con})$ in both the cases.
\item [\rm (v)] In general, ${\rm Aut}(\V_{con})$ need not be Abelian even if $\bS=\{id\}$   {(\cite[Corollary 6.6]{Gai})}. The following example shows that even ${\rm Aut}(\delta)$ need not be Abelian in such a situation.
\end{enumerate}
}
\end{rem}
\begin{ex}\label{danst3}
{\em Let $\V_{con}$ be such that $k_i \neq k_j$ if $i \neq j$ for all $i,j \in \{2,3,\dots ,m\}$. Then ${\bS}_{\delta}=\{id\}$. Let $\delta = D_{con}$. Define $\theta , \alpha \in {\rm Aut}(\V_{con})$ as follows:
$$\theta (y_1,y_2, \dots ,y_m,z)=(y_1,{(-1)}^{k_m}y_2, y_3, \dots, y_{m-1}, {(-1)}^{k_2}y_m,z)$$ and
$$\alpha = {\rm Exp}\big{(}{(-1)}^{k_m-1}y_2D_{con}\big{)}.$$
By Proposition \ref{ptorus1}, $\theta \in {\bT}_{\delta}$. Clearly, $\alpha \in {\rm Aut}(\delta)$. Now,
$$\theta \alpha(z) = \theta (z+{(-1)}^{k_m-1}y_2Dz) = z+{(-1)}^{2k_m-1}y_2^{k_2+1}y_3^{k_3} \dots y_m^{k_m}$$
and
$$ \alpha \theta (z) = \alpha (z) = z+{(-1)}^{k_m-1}y_2Dz=z+{(-1)}^{k_m-1}y_2^{k_2+1}y_3^{k_3} \dots y_m^{k_m}.$$
If $k_m$ is odd, then $\theta \alpha \neq \alpha \theta$. Hence ${\rm Aut}(\delta)$ is non-Abelian.
}

\end{ex}

\section{Danielewski varieties}\label{iso_Dv}
Let $\V$ be a Danielewski variety (as defined in Definition \ref{danvar}) with
$$\bK[\V]=\dfrac{\bK[Y_1,Y_2, \dots ,Y_m,Z]}{(Y_1Y_2^{k_2} \dots Y_m^{k_m}-P(Y_2, \dots, Y_m,Z)},$$
 such that $k_i (\in \bN)\geqslant 2$ for all $i \in \{2, \dots ,m\}$. Let $y_1,y_2, \dots y_m,z$ denote the images of $Y_1, \dots , Y_m, Z$ respectively in $\bK[\V]$. Furthermore, by a suitable change of coordinates, one may assume that the coefficient of $Z^{d-1}$ in $P$ is zero. By \cite[Lemma 6.2~(2)]{Gai}, $\V$ is irreducible, almost rigid and the {\it canonical lnd} $D_{\V}$ on $\V$ is given by
$$D_{\V}(y_1)=\frac{\partial P}{\partial z},~~D_{\V}(z)=y_2^{k_2}\dots y_m^{k_m}~~\text{and}~~D_{\V}(y_i)=0,~\forall i \geqslant 2.$$

\smallskip
\noindent
The symmetric group $S_m$ acts on $\bK[Y_2,\dots ,Y_m]$ by permutation of coordinates. The stabilizer $S(\V) \subseteq {\mathcal S}_m$ of the monomial $Y_2^{k_2}\dots Y_m^{k_m}$ permutes the $Y_i$'s with the same $k_i$'s. There is an action of the semi-direct product $G=S(\V) \rightthreetimes ({\bK^*})^{m}$ on \\$\bK(y_2,\dots ,y_m,z)(\cong {\bK}^{(m)})$ given by 
$$(\sigma, t_2, \dots t_{m+1}) \cdot (y_2, \dots y_m, z)=(t_2y_{\sigma (2)}, \dots ,t_my_{\sigma (m)}, t_{m+1}z).$$
By mapping $y_1 \rightarrow \dfrac{P(y_2,\dots ,y_m,z)}{y_2^{k_2} \dots y_m^{k_m}}$, one has a natural embedding $\varphi : \bK[\V] \rightarrow \bK(y_2, \dots,y_m,z)$. Let $\Ga~(\subseteq G)$ be the stabilizer of $\varphi(\bK[\V])$. Then one has a $\Ga$-action on $\V$. As shown in \cite[Section 7]{Gai} one should note that for any $g=(\sigma, t_2, \dots ,t_{m+1}) \in \Ga$, one can define $g \cdot y_1$ correctly, i.e., $g \cdot y_1 \in \bK[y_2,\dots y_m,z]$ if and only if 

\begin{equation}
 {\underset{j=2}{\overset{m}\prod}}y_j^{k_j} \mid s_i(t_2y_{\sigma (2)}, \dots ,t_m{\sigma (m)})-t^{d-i}_{m+1}s_i(y_2, \dots y_m)~~\forall i \in \{0,1,\dots ,d-2\}
\end{equation}
\smallskip
\noindent
$\Ga$ is called the {\it canonical group} of $\V$. Let us also denote the subgroup of ${\rm Aut}(\V)$ generated by the automorphisms induced by the canonical group by $\Ga$. The following theorem \cite[Theorem 7.11]{Gai} describes the structure of ${\rm Aut}(\V)$. 

\begin{thm}\label{danv1}
Let $\V$ be a Dainelewski variety. Then ${\rm Aut}(\V) \cong \Ga \rightthreetimes \U(D_{\V})$, where $\Ga$ is the canonical group of $\V$ and $\U(D_{\V})$ is the big unipotent group corresponding to the canonical {\it lnd} on $\V$.
\end{thm}

\smallskip
\noindent
As a result of Theorem \ref{danv1}, we have the following corollary.
\begin{cor}\label{danv3}
Let $\V$ be a Danielewski variety. For $\delta =hD_{\V}$, where $h \in \bK[y_2, \dots ,y_m]$, we have 
$${\rm Aut}(\delta) \cong {\Ga}_{\delta} \rightthreetimes \U(D_{\V}).$$
\end{cor}

\begin{proof} The result follows from Theorem \ref{danv1}.
\end{proof}
The next proposition is similar to Proposition \ref{ptorus1} of the previous section.
\begin{prop}\label{danv2}
Let $\delta =hD_{\V},$ where $h \in \bK[y_2,\dots ,y_m]$ and $\theta \in \Ga $ be given by $$\theta (y_2, \dots y_m,z)=(t_2y_{\sigma (2)}, \dots ,t_my_{\sigma (m)}, t_{m+1}z),$$ where $\sigma \in S(\V)$ and $(t_2, \dots ,t_{m+1}) \in {(\bK^*)}^{m}$. Then the following statements are equivalent.
\begin{enumerate}
\item [\rm (a)] $\delta \theta (y_1)=\theta \delta (y_1)$.
\item [\rm (b)] $t_{m+1}h=\theta(h){\overset{m}{\underset{j=2}\prod}} t_j^{k_j}$.
\item [\rm (c)] $\delta \theta (z)=\theta \delta (z)$.
\item [\rm (d)] $\theta \in {\rm Aut}(\delta)$.
\end{enumerate}
Moreover, if $\theta(h)=h$ $($in particular, when $h \in {\bK}^*)$, the statement ${\rm(b)}$ is equivalent to \\
${\rm(b)}^{\prime} ~ t_{m+1}={\underset{j=2}{\overset{m}\prod}} t_j^{k_j}$.\\
\end{prop}

\begin{proof}

For convenience, set ${\underset{j=2}{\overset{m}\prod}}t_j^{k_j}=\tilde{t}$ and ${\underset{j=2}{\overset{m}\prod}}y_j^{k_j}=\tilde{y}$. Then $\theta (\tilde{y})={\tilde{t}}{\tilde{y}}$. 

\noindent
{\bf \underline{(a)$\Leftrightarrow$(b)}} We observe that
$$\tilde{t}\tilde{y} \delta \theta (y_1)=\delta \theta (y_1\tilde{y})=\delta \theta \big{(}P(y_2, \dots ,y_m,z)\big{)}=\theta (\frac{\partial P}{\partial z})\delta (t_{m+1}z)=\theta (\frac{\partial P}{\partial z})t_{m+1}h\tilde{y}$$ and
$$\tilde{t}\tilde{y}\theta\delta(y_1)=\tilde{t}\tilde{y}\theta (h\frac{\partial P}{\partial z})=
\theta (\frac{\partial P}{\partial z})\theta(h)\tilde{t}\tilde{y}.$$
Now, (a)$\Leftrightarrow$(b) follows from the fact that $\bK[\V]$ is a domain.

\smallskip
\noindent
{\bf \underline{(b)$\Leftrightarrow$(c)}} Suppose (b) holds. Then
$$\delta \theta (z)=\delta (t_{m+1}z)=(t_{m+1}h)\tilde{y}=\big{(}\theta(h)\tilde{t}\big{)}\tilde{y}=\theta(h\tilde{y})=\theta \delta (z).$$ So $(b) \Rightarrow (c)$. Conversely, suppose (c) holds. Then,
$$\tilde{t}\tilde{y}\theta(h)=\theta(h\tilde{y})=\theta(\delta(z))=\delta \theta(z)=t_{m+1}h\tilde{y}.$$ Since $\bK[V]$ is a domain, we have $(c) \Rightarrow (b)$.

\smallskip
\noindent
{\bf \underline{(d)$\Leftrightarrow$(a)}}
The proof of (d)$\Rightarrow$(a) is obvious. Suppose (a) holds. Then, for any $j \in \{2,\dots ,m\}$ we have
$$ \delta \theta (y_j) = \delta(t_jy_{\sigma (j)})=0=\theta\delta (y_j).$$  Since $(a) \Leftrightarrow (c)$, we have $\theta \in {\rm Aut}(\delta)$, proving $(a) \Rightarrow (d)$.

\smallskip
\noindent
 If $\theta(h)=h$, then the equivalence of (b) and (b)$^{\prime}$ is obvious.
\end{proof}

\begin{rem}{\em
 If $m=2$, then the Danielewski variety turns out to be the generalised Danielewski surface discussed in Section \ref{iso_GDs}, provided $P$ breaks up into distinct linear factors.}
 \end{rem}
 
 \begin{rem}\label{danrem}
 {\em It may be possible, that $\Ga={\Ga}_{\delta}=\{id\}$ (for example, when $\bK[\V]=\frac{\bK[Y_1,Y_2,Z]}{(Y_1Y_2^2-(Z^3+(Y_2+1)Z+1))}$ \cite[Example 7.6]{Gai}). In such cases, ${\rm Aut}(\delta)=\U(D_{\V}).$
 }
 \end{rem}

 \section{A threefold given by Finston and Maubach}\label{iso_FM}
 Let $R=\dfrac{\bC[X,Y,Z]}{(X^a+Y^b+Z^c)}$, where $a,b,c$ are three pairwise coprime positive integers with $1/a+1/b+1/c<1$ and set $B_{m,n}=\dfrac{R[U,V]}{(X^mU-Y^nV-1)}$, where $m,n \in \bN$, $m,n \geqslant 2$.
 
\medskip
\noindent 
The above class of threefolds were introduced by D. Finston and S. Maubach in \cite{FM}. It was shown that $B_{m,n}$ is a UFD (\cite[Proposition 1]{FM}) and $B_{m,n}$ is an almost rigid domain with Makar-Limanov invariant ${\rm ML}(B_{m,n})=R$. Let $x,y,z,u,v$ be the images of $X,Y,Z,U,V$ respectively in $B_{m,n}$. The {\it canonical lnd} $D_{m,n}$ is given by (see \cite[Lemma 4]{FM})$$D_{m,n}(u)=y^n~~\text{and}~~D_{m,n}(v)=x^m.$$ The structure of ${\rm Aut}(B_{m,n})$ has also been described (\cite[Theorem 3]{FM}). We state the precise result below.

\begin{thm}\label{FM_Aut}
	${\rm Aut}(B_{m,n})$ is generated by the following automorphisms.
\begin{enumerate}

\item[\rm(1)] $\theta_f^{+}(x,y,z,u,v)={(x,y,z,u+fy^n,v+fx^m)}$ for $f\in R$.
\item[\rm(2)]$\theta_{\mu}^{*}(x, y, z, u, v)={({\mu}^ {bc}x, {\mu}^{ac}y, {\mu}^{ab}z, {\mu}^{-mbc}u, {\mu}^{-nac}v)}$ for $\mu\in{\bC^*}$.
\end{enumerate}
Moreover, ${\rm Aut}(B_{m,n})=\bC^*\ltimes (R, +)$. 
\end{thm}

\begin{rem}{\em Following the notation of Theorem \ref{FM_Aut}, we can observe that
\begin{enumerate}
\item[\rm(i)] for each $f\in{R}$, the automorphisms $\theta_f^{+}={\rm Exp}(fD_{m,n})~(\in{\U(D_{m,n})})$ and hence is in ${\rm Aut}(fD_{m,n})$, for all $f\in R$,
\item[\rm(ii)]for each $\mu\in{\bC^*}$, $\theta_{\mu}^{*}$ corresponds to an action of the one-dimensional algebraic torus $\bC^*$ on $B_{m,n}$. 
\end{enumerate}  
}
\end{rem}

\noindent
Let $\bT$ denote the set of all automorphisms induced by the action of $\bC^{*}$ on $B_{m,n}$. The next lemma gives necessary and sufficient conditions for an automorphism of $B_{m,n}$ to be an element of $\bT_{\delta}$ for any {\it lnd} $\delta$ on $B_{m,n}$.

\begin{lem}\label{FM_isoaut}
	Let $\delta=hD_{m,n}$, where $h\in R=\bC[x,y,z]$ and $\theta:=\theta_{\mu}^{*}\in{\bT}$. Then the following statements are equivalent.
	\begin{enumerate}
		\item [\rm(i)] $\theta\in{\rm Aut}(\delta)$.
		\item [\rm(ii)] $\delta\theta(u)=\theta\delta(u)$.
		\item [\rm(iii)] $\delta\theta(v)=\theta\delta(v)$.
		\item [\rm(iv)] ${\mu}^{nac+mbc}\theta(h)=h$.
	\end{enumerate}
\end{lem}
\begin{proof}

\medskip
\noindent
	{\bf \underline {(ii)$\Leftrightarrow$(iv)}}
	We have 
	
	\smallskip
	\noindent
	$$\delta\theta(u)=\delta ({\mu}^{-mbc}u)={\mu}^{-mbc}hy^n$$ and
	$$\theta\delta(u)=\theta(hy^n)=\theta(h){\mu}^{nac}y^n.$$
	Since $B_{m,n}$ is a domain, (ii) holds if and only if (iv) holds.

\smallskip
\noindent
{\bf \underline{(iii)$\Leftrightarrow$(iv)}}
    Note that
	
	\smallskip
	\noindent
	$$\delta\theta(v)=\delta ({\mu}^{-nac}v)={\mu}^{-nac}hx^m$$ and
	$$\theta\delta(v)=\theta(hx^m)=\theta(h){\mu}^{mbc}x^m.$$
	Since $B_{m,n}$ is a domain, (iii) holds if and only if (iv) holds.
	
\smallskip
\noindent
{\bf \underline {(i)$\Leftrightarrow$(ii)}}

\smallskip
\noindent
Clearly, (i)$\Rightarrow$(ii). Since (ii)$\Rightarrow$(iii) and $x,y,z\in{\rm Ker}(\delta)$, hence (ii)$\Rightarrow$(i).   

\end{proof} 

We finish this section by describing the structure of the group ${\rm Aut}(\delta)$.
\begin{thm}\label{thmfm}
	Let $\delta=hD_{m,n}$, where $h\in R$. Then, ${\rm Aut}(\delta)\cong{G}\ltimes{(R,+)}$, where $G$ is a finite cyclic group. 
\end{thm}
\begin{proof}
	By Theorem \ref{FM_Aut}, ${\rm Aut}(\delta)\cong{\bT_{\delta}\ltimes (R,+)}$.
	 We claim that $\bT_{\delta}$ is a finite cyclic group. Let $h={\underset {r,s,t\geqslant 0, r<a }\sum}a_{r,s,t}x^ry^sz^t$. Then the condition (iv) in Lemma \ref{FM_isoaut}, will imply $\theta\in{\rm Aut}(\delta)$ if and only if $\mu^{bc(m+r)+ac(n+s)+abt}=1$, for each $(r,s,t)$ occurring in the expression of $h$. 
\end{proof}

\begin{rem}\label{fmrem}{\em Following the notations of Theorem \ref{thmfm}, we make the following observations.
\begin{enumerate}
\item [\rm (i)] If $\delta=D_{m,n}$, then $\U(\delta)$ is properly contained in ${\rm Aut}(\delta)$ since $mbc+nac > 2$.
\item [\rm (ii)] Since $\frac{1}{a}+\frac{1}{b}+\frac{1}{c} < 1$, we may assume $a > 2$. Let $\delta=(x+x^2)D_{m,n}$. Then $\theta_{\mu}^* \in {\bT}_{\delta}$ if and only if $\mu=1$. So ${\rm Aut}(\delta)=\U(\delta)$.
\end{enumerate}
}
\end{rem}

\section{Double Danielewski surfaces}\label{iso_DDs}

In this section we consider the Double Danielewski surfaces introduced by N. Gupta and S. Sen in \cite{GS1}. We give necessary and sufficient conditions for an automorphism of these surfaces to be an element of an isotropy subgroup. We first recall some definitions and fix some notations which will be used throughout this section. 

\begin{defn}\label{DDs}
{\em Let $B$ be a $k$-algebra. $B$ is said to be a {\it Double Danielewski surface over $k$} if $B$ is isomorphic to the $k$-algebra $$B_{d_1,P_1;d_2,P_2}:=\dfrac{k[X,Y_1,Y_2,Y_3]}{\big{(}X^{d_1}Y_2-P_1(X,Y_1), X^{d_2}Y_3-P_2(X,Y_1,Y_2)\big{)}},$$ where for each $i=1,2$, $d_i\geqslant{2}$, $P_i$ is monic in $Y_i$ and $r_i:={\rm deg}_{Y_i}(P_i)\geqslant{2}$.} 
\end{defn}

\smallskip
\noindent
By $x,y_1,y_2,y_3$ we will denote the images of $X,Y_1,Y_2,Y_3$ respectively in $B_{d_1,P_1;d_2,P_2}$. In \cite[Theorem 2.1]{GS2}, N. Gupta and S. Sen proved that $B_{d_1,P_1;d_2,P_2}$ is an almost rigid domain with the {\it canonical lnd} $$D_{d_1,P_1;d_2,P_2}=x^{d_1+d_2}\dfrac{\partial}{\partial y_1}+x^{d_2}\dfrac{\partial P_1}{\partial y_1}\dfrac{\partial}{\partial y_2}+(\dfrac{\partial P_2}{\partial y_1}x^{d_1}+\dfrac{\partial P_2}{\partial y_2}\dfrac{\partial P_1}{\partial y_1}) \dfrac{\partial}{\partial y_3}.$$

\smallskip 
The following result, proved by N. Gupta and S. Sen in \cite[Theorems 3.13, 3.12]{GS1}, gives a necessary and sufficient condition for an $k$-algebra endomorphism of $B_{d_1,P_1;d_2,P_2}$ to be an automorphism.
\begin{thm}\label{Aut_DD}
 Let $\alpha:B_{d_1,P_1;d_2,P_2}\longrightarrow{B_{d_1,P_1;d_2,P_2}}$ be a $k$-algebra endomorphism. Then $\alpha\in{\rm Aut}(B_{d_1,P_1;d_2,P_2})$ if and only if the following conditions hold.
 \begin{enumerate}
  \item[\rm(i)]$\alpha(x)=\lambda x$, for some $\lambda\in{k^*}$.
  \item[\rm(ii)]$\alpha(y_1)=ay_1+ b(x)$, for some $a\in{k^*}$ and $b(x)\in{k[x]}$.
 \end{enumerate}
\end{thm}
Now, we will prove the main result of this section in which we give necessary and sufficient conditions for an automorphism of $B_{d_1,P_1,d_2,P_2}$ to be an element of the isotropy subgroup of any replica of the canonical {\it lnd}.
\begin{prop}\label{y_1}
Let $\delta:=f(x)D_{d_1,P_1;d_2,P_2}$, where $f(x)={\underset{i=0}{\overset{l}\sum}}{a_ix^{n_i}}\in{k[x]}~(\text{for each} ~i,~a_i \in k^*$, $l,n_i \in \bN \cup \{0\})$. Suppose $\alpha\in{\rm Aut}(B_{d_1,P_1;d_2,P_2})$ with $\alpha(x)=\lambda x$ and $\alpha(y_1)=ay_1+b(x)$, for some $\lambda,a\in{k^*}$, $b(x)\in{k[x]}$. Then the following are equivalent.
\begin{enumerate}
\item[\rm(i)] $\alpha\in {\rm Aut}(\delta)$.
\item[\rm(ii)]$\delta\alpha(y_1)=\alpha \delta(y_1)$.
\item[\rm(iii)]$\delta\alpha(y_2)=\alpha \delta(y_2)$.
\item[\rm(iv)] $a=\lambda^{d_1+d_2+n_{0}}=\lambda^{d_1+d_2+n_{i}}$, for all $i$, $1\leqslant i\leqslant l$.
\end{enumerate}
In particular, if $\delta=D_{d_1,P_1;d_2,P_2}$, then $\alpha\in{\rm Aut}(\delta)$ if and only if $a=\lambda^{d_1+d_2}$.

\end{prop}
\begin{proof}

\noindent 
Since $x^{d_1}y_2=P(x,y_1)$, $x^{d_2}y_3=Q(x,y_1,y_2)$, we have
\begin{enumerate}
\item[\rm(1)] $\lambda^{d_1}x^{d_1}\alpha(y_2)=P_1\big{(}\alpha(x),\alpha(y_1)\big{)}$,
\item[\rm(2)] $\lambda^{d_2}x^{d_2}\alpha(y_3)=P_2\big{(}\alpha(x),\alpha(y_1),\alpha(y_2)\big{)}$,
\item[\rm(3)] $x^{d_1}\delta(y_2)=\dfrac{\partial {P_1}}{\partial {y_1}}\delta(y_1)$ and
\item[\rm(4)] $x^{d_2}\delta(y_3)=\dfrac{\partial {P_2}}{\partial{y_1}}\delta(y_1)+\dfrac{\partial {P_2}}{\partial{y_2}}\delta(y_2)$.
\end{enumerate}

\noindent
Since $\rm(i)\Rightarrow \rm(ii)$, $\rm(i)\Rightarrow \rm(iii)$ are clear, it is enough to show that $\rm(ii)\Leftrightarrow \rm(iii)$, $\rm(ii)\Leftrightarrow \rm(iv)$ and $\rm(ii)\Rightarrow \rm(i)$.

\medskip
\noindent
{\bf \underline{(ii)$\Leftrightarrow$(iii)}}
We have

\smallskip
\noindent
$\delta\alpha(y_2)=\alpha \delta(y_2)\Leftrightarrow{\lambda^{d_1}x^{d_1}\alpha \delta(y_2)=\lambda^{d_1}x^{d_1}\delta\alpha(y_2)}\Leftrightarrow{\alpha\rbo x^{d_1}\delta(y_2)\rbc=\delta\rbo\lambda^{d_1}x^{d_1}\alpha(y_2)\rbc}$. Hence by (1) and (3), $\delta\alpha(y_2)=\alpha \delta(y_2)\Leftrightarrow{\alpha\rbo\dfrac{\partial P_1}{\partial {y_1}}\delta(y_1)\rbc=\delta\Big{(}P_1\big{(}\alpha(x),\alpha(y_1)\big{)}\Big{)}}$. Therefore by Lemma \ref{chainrule}, $\delta\alpha(y_2)=\alpha \delta(y_2)\Leftrightarrow{\alpha(\dfrac{\partial P_1}{\partial {y_1}})\alpha \delta(y_1)=\alpha(\dfrac{\partial P_1}{\partial {y_1}})\delta\alpha({y_1})}$.

\smallskip
\noindent
Since $\dfrac{\partial P_1}{\partial {y_1}}\neq{0}$, $\delta\alpha(y_2)=\alpha \delta(y_2)\Leftrightarrow{\delta\alpha(y_1)=\alpha \delta(y_1)}$. 

\medskip
\noindent
{\bf \underline{(ii)$\Rightarrow$(i)}}

\smallskip
\noindent
Since $\delta\alpha(x)=0=\alpha \delta(x)$, it is enough to show that $\delta\alpha(y_3)=\alpha \delta(y_3)$.

\smallskip
\noindent
Considering the images under $\delta$ of both sides of (2), we obtain
\begin{align*}
\lambda^{d_2}x^{d_2}\delta\alpha(y_3)&=\alpha(\dfrac{\partial P_2}{\partial {y_1}})\delta\alpha(y_1)+\alpha(\dfrac{\partial P_2}{\partial {y_2}})\delta\alpha(y_2)\\
&=\alpha(\dfrac{\partial P_2}{\partial {y_1}})\alpha \delta(y_1)+\alpha(\dfrac{\partial P_2}{\partial {y_2}})\alpha \delta(y_2)\\
&=\alpha\big{(}\dfrac{\partial P_2}{\partial {y_1}}\delta(y_1)+\dfrac{\partial P_2}{\partial {y_2}}\delta(y_2)\big{)}\\&=\alpha\rbo x^{d_2}\delta(y_3)\rbc,~ \text{by (4)}\\
&=\lambda^{d_2}x^{d_2}\alpha \delta(y_3).
\end{align*}
Hence, $\delta\alpha(y_3)=\alpha \delta(y_3)$.

\medskip
\noindent
{\bf \underline{(ii)$\Leftrightarrow$(iv)}}
We have

\smallskip
\noindent
$\delta\alpha(y_1)=af(x)x^{d_1+d_2}$ and $\alpha\delta (y_1)=\lambda^{d_1+d_2}f(\lambda x)x^{d_1+d_2}$. Now, by Lemma \ref{root of unity}, $\delta\alpha(y_1)=\alpha \delta(y_1)$ if and only if $a=\lambda^{d_1+d_2+n_{i}}$, for all $i$, $0\leqslant i\leqslant l$.
\end{proof} 

\begin{rem} {\em From Proposition \ref{y_1}, we can make the following observations. 
\begin{enumerate}
\item [\rm (i)] Without loss of generality, we can assume that $n_0 <n_1 < \dots < n_l$. If $n_0=0$ (equivalently, $f(0) \neq 0$), then by (iv) Proposition \ref{y_1}, $\lambda$ is a $p$-th root of unity, where $p:={\rm GCD}(n_1, \dots ,n_l)$. Similarly, if $n_0 > 0$ (equivalently, $f(0)=0$), then $\lambda$ is a $q$-th root of unity, where $q:={\rm GCD}(\{n_{i+1}-n_i~:~0 \leqslant i \leqslant l-1\})$.
\item [\rm (ii)] When $\delta=D_{d_1,P_1;d_2,P_2}$, any $\lambda~ (\neq 1) \in k^*$ gives an element of ${\rm Aut}(\delta)$ (such that $\alpha(x)=\lambda x$) which is not in $\U(\delta)$. So in this case, ${\rm Aut}(\delta)$ properly contains $\U(\delta)$. However, this may not be true when $\delta$ is a replica of $D_{d_1,P_1;d_2,P_2}$. For example, if $f(x)=x+x^2$, then $\lambda=a=1$. So, if $\alpha \in {\rm Aut}(\delta)$, then $\alpha \in \U(\delta)$. Hence ${\rm Aut}(\delta)=\U(\delta)$.
\end{enumerate}
}
\end{rem}

\noindent
\begin{center}
{\bf Acknowledgement:}
\end{center}

\noindent
The authors thank Professor Ivan Arzhantsev for his valuable comments while going through the earlier drafts and suggesting 
improvements.

 	{\small{
 
}}

\end{document}